\RequirePackage{fix-cm}
\documentclass[11pt]{article} 
\usepackage{graphicx}
\usepackage{float,pict2e}
\usepackage[<pict2e options>]{curve2e}
\usepackage{amsmath}
\usepackage{amssymb}

\newtheorem{theorem}{\bf Theorem}[section]
\newtheorem{lemma}{\bf Lemma}[section]

\begin{document}

\title{On triangle-free graphs of order 10 with prescribed 1-defective chromatic number 
}


\author{Nirmala Achuthan \and N.R. Achuthan \and G. Keady.\\
Department of Mathematics and Statistics, Curtin University\\
Bentley  WA 6102,  \textsc{Australia}
}



\date{\today} 

\maketitle

\begin{abstract}
A graph is $(m,k)$-colourable if its vertices can be coloured with
$m$ colours such that the maximum degree of any subgraph induced on
vertices receiving the same colour is at most $k$.
The $k$-defective chromatic number for a graph 
is the least positive integer $m$ for which the graph  is
$(m,k)$-colourable. All triangle-free graphs on 8 or fewer vertices
are $(2,1)$-colourable. There are exactly four triangle-free graphs
of order 9 which have 1-defective chromatic number 3. We show that
these four graphs appear as subgraphs in almost all triangle-free
graphs of order 10 with 1-defective chromatic number equal to 3. In
fact there is a unique triangle-free $(3,1)$-critical graph on 10
vertices and we exhibit this graph. 
\end{abstract}


\par\noindent{\bf  Math Review Codes.}
{MSC 05C15 \and MSC 05C35}

\section{Introduction}
\label{intro}
We consider in this paper undirected graphs with no
loops or multiple edges. For all undefined concepts and terminology
we refer to~\cite{CLZ}.

Given a graph $G$, $d_G(u)$, $N_G(u)$ and $N_G[u]$ denote
respectively the degree, the neighbourhood, and the closed
neighbourhood of a vertex $u$ in $G$. The union of graphs $G_1$ and
$G_2$ is denoted by $G_1\cup{G_2}$. For convenience we write $2G$ in
place of $G\cup{G}$.

Let $k$ be a nonnegative integer. A subset $U$ of the vertex set
$V(G)$ is {\it $k$-independent} if $\Delta(G[U])\le{k}$. A
$0$-independent set is an independent set in the usual sense. A
graph $G$ is {\it $(m,k)$-colourable} if it is possible to assign
$m$ colours, say $1,2,\ldots,m$ to the vertices of $G$, one colour
to each vertex, such that the set of all vertices receiving the same
colour is $k$-independent. The smallest integer $m$ for which $G$ is
$(m,k)$-colourable is called the {\it $k$-defective chromatic
number} of $G$ and is denoted by $\chi_k(G)$. A graph $G$ is said to
be {\it $(m,k)$-critical} if $\chi_k(G)=m$ and $\chi_k(G-u) < m$ for
every $u$ in $V(G)$. A graph $G$ is said to be {\it
$(m,k)$-edge-critical} if $\chi_k(G)=m$ and $\chi_k(G-e) < m$ for
every $e$ in $E(G)$.

It is easy to see that the following statements are equivalent.
\begin{enumerate}
\item[(i)] $G$ is $(m,k)$-colourable.
\item[(ii)] There exists a partition of $V(G)$ into $m$ sets each of which is
$k$-independent.
\item[(iii)] $\chi_k(G)\le{m}$.
\end{enumerate}
Note that $\chi_0(G)$ is the usual chromatic number. It is easy to
see that $\chi_k(G)\le\lceil\frac{|V(G)|}{k+1}\rceil$. The concept
of $k$-defective chromatic number has been extensively studied in
the literature (see~\cite{AAS11,CGJ,Fr,GH,HS,SAA97b,Wo}). Given a
positive integer $m$, it is well known that there exists a
triangle-free graph $G$ with $\chi_k(G)=m$. A natural question that
arises is: what is the smallest order of a triangle-free graph $G$
with $\chi_k(G)=m$? We denote this smallest order by $f(m,k)$. The
parameter $f(m,0)$ has been studied by several authors
(see~\cite{Avis,Chvatal,JR,HM}) and  $f(m,0)$ is determined for
$m\le{5}$. It has also been shown that $f(3,1)=9$ and $f(3,2)=13$.
Furthermore the corresponding extremal graphs have been
characterized (see ~\cite{SAA97b,AAS11}).

In this paper we characterize triangle-free graphs of order 10 with
$\chi_1(G)=3$.
In a subsequent paper~\cite{AAK12b}  
we build from the results of this paper to
determine the smallest order of a triangle-free  planar  graph which has 
1-defective chromatic number 3.

In all the figures in this paper a double line between sets $X$ and
$Y$ means that every vertex of $X$ is adjacent to every vertex of
$Y$.

\section{Preliminary results}
\label{sec:1}

We need the following results, proofs of the theorems being in the
papers cited.

\begin{theorem} (\cite{HS,Lovasz})
\label{thm:Lovasz} Let $G$ be a graph with maximum degree $\Delta$.
Then
$$ \chi_k(G) \le \lceil \frac{\Delta+1}{k+1}\rceil
= 1 +  \lfloor\frac{\Delta}{k+1}\rfloor .
$$
\end{theorem}

\smallskip

\begin{figure}[H]
 \centering
 \setlength{\unitlength}{0.7cm}
 \begin{picture}(15,10)(0,-7.3)
 \multiput(0,0)(8,0){2}{%
  \multiput(0,0)(0,-5){2}{%
  \put(0,0){\circle{0.7}
            \makebox(0,0)[r]{$u$}}
  \multiput(1,1)(4,0){2}{%
       \polyline(0,0)(2,0)(2,1)(0,1)(0,0)
       \multiput(0.5,0.5)(1,0){2}{\circle*{0.12}}}
  \put(1,1){%
       \put(0,-3){%
            \multiput(0.5,0.5)(1,0){2}{\circle*{0.12}}}}
  \put(3,1.5){%
            \multiput(0,0.05)(0,-0.1){2}{\line(1,0){2}}}
  \put(1,1){%
       \put(0.5,0.6){\makebox(0,0)[b]{$u_1$}}
       \put(1.5,0.6){\makebox(0,0)[b]{$u_2$}}}
  \put(5,1){%
       \put(0.5,0.6){\makebox(0,0)[b]{$z_1$}}
       \put(1.5,0.6){\makebox(0,0)[b]{$z_2$}}}
  \put(1,-2){%
       \put(0.5,0.4){\makebox(0,0)[t]{$u_3$}}
       \put(1.5,0.4){\makebox(0,0)[t]{$u_4$}}}
  \polyline(.240,.254)(1,1.394)
  \polyline(.142,.320)(1,1.606)
  \polyline(.240,-.254)(1,-1.394)
  \polyline(.142,-.320)(1,-1.606)}}

 \multiput(0,0)(8,0){2}{%
  \put(5,-2){%
            \multiput(0.5,0.5)(1,0){2}{\circle*{0.12}}}
  \multiput(1,1)(4,0){2}{%
       \put(0,-3){%
            \polyline(0,0)(2,0)(2,1)(0,1)(0,0)}}
  \put(3,-1.5){%
            \multiput(0,0.05)(0,-0.1){2}{\line(1,0){2}}}
  }
  \put(0,0){%
    \put(5,-2){%
       \put(0.4,0.4){\makebox(0,0)[rt]{$z$}}
       \put(1.5,0.4){\makebox(0,0)[t]{$z_3$}}}
    \polyline(5.5,1.5)(5.5,-1.5)(6.5,1.5)
    \put(3.5,-2.2){\makebox(0,0)[t]{$G_1$}}}
  \put(8,0){%
    \put(5,-2){%
       \put(0.5,0.4){\makebox(0,0)[t]{$z$}}
       \put(1.5,0.4){\makebox(0,0)[t]{$z_3$}}}
    \polyline(6.5,-1.5)(5.5,1.5)(5.5,-1.5)(6.5,1.5)
    \put(3.5,-2.2){\makebox(0,0)[t]{$G_2$}}}
  \put(0,-5){%
    \put(5,-2){%
            \multiput(0.5,0.5)(1,0){2}{\circle*{0.12}}}
    \multiput(1,1)(4,0){2}{%
       \put(0,-3){%
            \polyline(0,0)(2,0)(2,1)(0,1)(0,0)}}
    \put(3,-1.5){%
            \multiput(0,0.05)(0,-0.1){2}{\line(1,0){2}}}
    \put(5,-2){%
       \put(0.5,0.4){\makebox(0,0)[t]{$z$}}
       \put(1.5,0.4){\makebox(0,0)[t]{$z_3$}}}
       \multiput(5.95,1)(0.1,0){2}{\line(0,-1){2}}
    \put(3.5,-2.2){\makebox(0,0)[t]{$G_3$}}}
  \put(8,-5){%
    \put(1,-2){%
            \polyline(0,0)(3,0)(3,1)(0,1)(0,0)}
    \put(1,-2){%
       \put(2.5,0.5){\circle*{0.12}}
       \put(2.5,0.4){\makebox(0,0)[t]{$u_5$}}}
    \multiput(0,0)(0,0.1){2}{%
       \polyline(4,-1.55)(5.654,-1.55)}
    \multiput(0,0)(0.1,0){2}{%
       \polyline(5.950,-1.154)(5.950,1)}
    \put(6,-1.5){\circle{0.7}
            \makebox(0,0)[r]{$z$}}
    \put(3.5,-2.2){\makebox(0,0)[t]{$G_4$}}}
 \end{picture}
\caption{\label{fig:one} The critical graphs of order 9 with
$\chi_1(G)=3$: $G_1$ to $G_4$ of~\cite{SAA97b}. }

\end{figure}
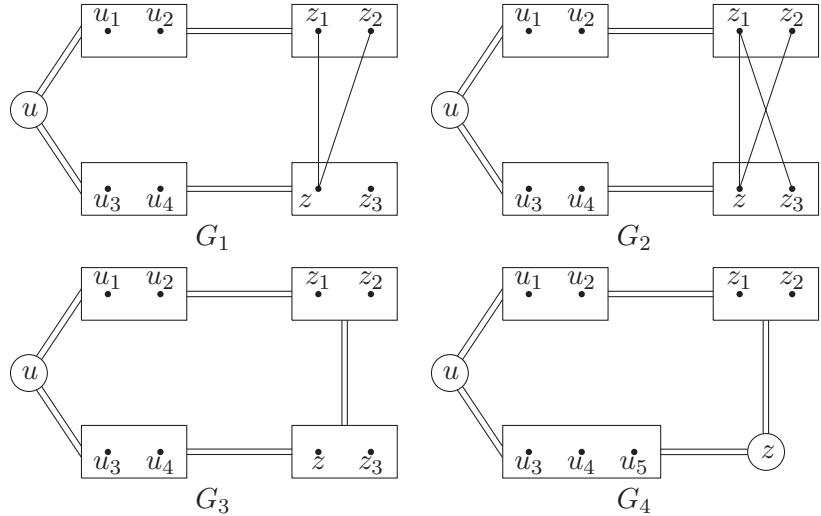

\begin{theorem} (\cite{SAA97b})
\label{thm:SAAb} The smallest order of a triangle-free graph with
$\chi_1(G)=3$ is 9, that is, $f(3,1)=9$. Moreover, $G$ is a
triangle-free graph of order 9 with $\chi_1(G)=3$ if and only if it
is isomorphic to one of the graphs $G_i$, $1\le{i}\le{4}$ given in
Figure \ref{fig:one}.
\end{theorem}

\section{Main results}
\label{sec:2}

Consider a graph $G$ of order $n$. The following notation is used
repeatedly in the paper:
\begin{eqnarray}
u\  {\mbox{\rm\  is a vertex of degree\ }} \Delta(G), &\qquad&
A=N_G(u),\qquad B=V(G)-N_G[u],\qquad
\label{eq:ABdef}\\
H=G[B] & {\mbox{\rm and\ }} & z\in{B} {\mbox{\rm\  with\ }}
d_H(z)=\Delta(H). \label{eq:Hdef}
\end{eqnarray}

We henceforth denote the vertex set $V(G)$ by $V$ and the edge set
$E(G)$ by $E$.
\medskip

\begin{lemma}
\label{lem:four} Let $G$ be a triangle-free graph. In the notation
described above, suppose that $\Delta(H)=|B|-1$ and
$|A\cap{N_G(z)|}\le{2 k}$, where $k$ is a nonnegative integer. Then
$\chi_k(G)\le{2}$.
\end{lemma}
\smallskip

\noindent{\it Proof.} Consider a partition of $A\cap{N_G(z)}$ into
two sets $A_{11}$ and $A_{12}$ such that $|A_{1i}|\le{k}$ for $i=1$
and $2$. Since $G$ is triangle-free, the sets
$N_H(z)\cup{\{u\}}\cup{A_{11}}$ and $(A-A_{11})\cup{\{z\}}$ are both
$k$-independent. Hence $\chi_k(G)\le{2}$. $\square$
\medskip

\begin{lemma}
\label{lem:five}
Let $G$ be a triangle-free 
graph of order 10 with $\chi_1(G)\ge{3}$. Then (i) $\Delta(H)\ge{2}$
and (ii) $4\le{\Delta(G)}\le{6}$.
\end{lemma}
\smallskip

\noindent{\it Proof.} The lower bound for $\Delta(G)$ follows from
Theorem~\ref{thm:Lovasz}. Let $u\in{V}$ with $d_G(u)=\Delta(G)$. If
$\Delta(H)\le{1}$, then $\{u\}\cup{B}$ is 1-independent. Since $A$
is also 1-independent, this implies $\chi_1(G)\le{2}$. Thus
$\Delta(H)\ge{2}$ and hence $|B|\ge{3}$ implying that
$\Delta(G)=|A|\le{6}$. $\square$
\medskip

\begin{lemma}
\label{lem:six} Let $G$ be a triangle-free graph of order 10 with
$\Delta(G)={6}$. If $\chi_1(G)={3}$ then there exists a vertex $u^*$
in $G$ such that $G-{u^*}\cong G_4$.
\end{lemma}
\smallskip

\noindent{\it Proof.} Assume that $\chi_1(G)={3}$. Using the
notation described before we have $|B|=3$. From (i) of
Lemma~\ref{lem:five} we have $\Delta(H)\ge{2}$. Thus $\Delta(H)=2$.

Let $z\in{B}$ with $d_H(z)=2$. Using Lemma~\ref{lem:four}, we
conclude that $|A\cap{N_G(z)}|\ge{3}$. Also, as $d_G(z)\le{6}$,
$|A\cap{N_G(z)}|\le{4}$.

Let $A_1=A\cap{N_G(z)}$, $A_2=A-A_1$ and $N_H(z)=\{z_1,z_2\}$. Since
$G$ is $K_3$-free, the set $A_1\cup{\{z_1,z_2\}}$ is 0-independent.
If $z_1$ is adjacent to at most one vertex of $A_2$, then
$$A\cup\{z_1\} {\mbox{\rm is 1-independent.\ \ So is\ \ }}
V-(A\cup\{z_1\})=\{u,z,z_2\}. $$ It follows that $\chi_1(G)\le{2}$,
a contradiction. Hence $z_1$ (similarly $z_2$) has at least two
neighbours in $A_2$. Since $|A_2|\le{3}$, $z_1$ and $z_2$ have at
least one common neighbour in $A_2$.

Suppose that there is exactly one common neighbour, say $x$, of
$z_1$ and $z_2$ in the set $A_2$. This implies that $|A_2|=3$ and
$X=(A-\{x\})\cup{\{z_1,z_2\}}$ is 1-independent. Since
$V-X=\{u,x,z\}$ is also 1-independent  we have $\chi_1(G)\le{2}$, a
contradiction. Thus $A_2$ has at least two common neighbours, say
$x$ and $y$, of $z_1$ and $z_2$.

Now select a vertex $u^*$ from $A$ as follows. If $|A_1|=4$ then
$u^*$ is any vertex of $A_1$. Otherwise, that is, if $|A_1|=3$ then
$u^*$ is a vertex in $A_2$ (note that $|A_2|=3$) different from $x$
and $y$. Now it is easy to verify that $G-u^*\cong{G_4}$. Hence the
result. $\square$

\medskip

\begin{lemma}
\label{lem:seven}
 Let $G$ be a triangle-free graph of order 10 with
${\Delta(G)}={5}$. If $\chi_1(G)={3}$ then either there exists a
vertex $u^*$ with $G-u^*\cong{G_i}$ for $1\le{i}\le{4}$ or
$G\cong{G_5}$ illustrated in Figure~\ref{fig:g5}.
\end{lemma}
\medskip

\begin{figure}[H]
 \centering
 \setlength{\unitlength}{0.7cm}
 \begin{picture}(8,5)(0.0,1.5)
  \put(0.5,3){%
  \polyline(0,0)(1.5,0)(1.5,1.5)(0,1.5)(0,0)}
 \put(0.7,4.2){\makebox(0,0)[l]{$u$}}
 \put(0.7,3.3){\makebox(0,0)[l]{$v$}}
 \put(1.2,4.2){\circle*{0.2}}
 \put(1.2,3.3){\circle*{0.2}}
  \put(3.5,3.5){%
  \polyline(0,0)(1.5,0)(1.5,3)(0,3)(0,0)}
\put(3.7,6.0){\makebox(0,0)[l]{$u_5$}}
 \put(3.7,5.0){\makebox(0,0)[l]{$u_4$}}
 \put(3.7,4.0){\makebox(0,0)[l]{$u_3$}}
 \put(4.3,6.0){\circle*{0.2}}
 \put(4.3,5.0){\circle*{0.2}}
 \put(4.3,4.0){\circle*{0.2}}
 \put(6.5,4.5){%
  \polyline(0,0)(1.5,0)(1.5,1.5)(0,1.5)(0,0)}
\put(7.6,5.6){\makebox(0,0)[r]{$z_2$}}
 \put(7.6,4.8){\makebox(0,0)[r]{$z_1$}}
 \put(7.1,5.6){\circle*{0.2}}
 \put(7.1,4.8){\circle*{0.2}}
\put(3.5,1.0){%
  \polyline(0,0)(1.5,0)(1.5,2)(0,2)(0,0)}
\put(3.7,1.5){\makebox(0,0)[l]{$u_1$}}
 \put(3.7,2.5){\makebox(0,0)[l]{$u_2$}}
 \put(4.2,1.5){\circle*{0.2}}
 \put(4.2,2.5){\circle*{0.2}}
 \put(7.25,2.0){\circle{0.8}
            \makebox(0,0)[r]{$z$}}
 \polyline(4.3,6.0)(7.1,5.6)(4.3,4.0)(7.1,4.8)(4.3,5.0)
 \polyline(2.0,4.23)(3.5,5.53)
  \polyline(2.0,4.05)(3.5,5.35)
  \polyline(2.0,3.7)(3.5,2.2)
  \polyline(2.0,3.5)(3.5,2.0)
 \polyline(7.25,4.5)(7.25,2.4)
  \polyline(7.1,4.5)(7.1,2.36)
  \polyline(5,2.15)(6.85,2.15)
  \polyline(5,2.0)(6.84,2.0)
 \end{picture}
  \vspace{0.3 cm}
\caption{\label{fig:g5} $G_5$}
\end{figure}

\noindent{\it Proof.} Suppose that $\chi_1(G)=3$. Using the notation
described before, it follows that $|B|=4$. Now using
Lemma~\ref{lem:four} and Lemma~\ref{lem:five}(i), we have
$\Delta(H)=2$.
Let $v\in{B}$ such that $(z,v)\not\in{E}$, $N_H(z)=\{z_1,z_2\}$ and
$A_1=A\cap{N_G(z)}$. Note that $|A_1|\le{3}$.

\noindent{\bf Case i.} $|A_1|=3$.\\
Let $A-A_1=\{x_1,x_2\}$. Suppose that
$(z_1,x_1)\not\in{E}$.\\
{\bf Claim~\ref{lem:seven}.1.} $(v,z_2)\in{E}$.\\
 Since $\chi_1(G)=3$ and
$$A\cup{\{z_1\}} {\mbox{\rm \   is 1-independent}},
V-(A\cup\{z_1\})= \{u,v,z,z_2\} {\mbox{\rm \   is not
1-independent}}.$$
This proves Claim~\ref{lem:seven}.1.\\
{\bf Claim~\ref{lem:seven}.2.} $(v,x_2)\in{E}$.\\
 Since $\chi_1(G)=3$ and
$(A-\{x_2\})\cup\{z_1,z_2\}$ is 1-independent, it follows that
$\{u,z,v,x_2\}$ is not 1-independent. This in turn implies that
$(v,x_2)\in{E}$.

\smallskip

Combining Claims~\ref{lem:seven}.1 and~\ref{lem:seven}.2 with the
assumption that $G$ is triangle-free, we have $(z_2,x_2)\not\in{E}$.
Now, note that the sets
$$X_1=A\cup\{z_1,z_2\} {\mbox{\rm \ and\ \ }}
V-X_1=\{u,z,v\} {\mbox{\rm \  are both 1-independent}}$$ implying
that $\chi_1(G)\le{2}$, a contradiction. Thus $(z_1,x_1)\in{E}$.
Using similar arguments we conclude that $(z_1,x_2)\in{E}$ and
$(z_2,x_i)\in{E}$ for $i=1$, 2. Now, clearly, $G-v\cong{G_4}$. This
completes Case i.
\smallskip

\noindent{\bf Case ii.} $|A_1|\le{2}$.\\
Since $\Delta(H)=2$ and $|B|=4$, clearly $H$ is either
$P_3\cup{K_1}$ or $P_4$ or $C_4$.

Let us first consider the case that $H\cong{P_3\cup{K_1}}$ or $P_4$.\\
If $|A_1|\le{1}$ then the sets $X=A\cup\{z\}$ and $V-X$ partition
the vertex set $V$ of $G$ into two 1-independent sets
implying that $\chi_1(G)\le{2}$, a contradiction.\\
Hence $|A_1|=2$. Let $A_1=\{u_1,u_2\}$.
If $(v,u_1)\not\in{E}$ then the sets $X_1=\{u,u_1\}\cup(B-\{z\})$
and $V-X_1$ partition $V$ into 1-independent sets. This implies that
$\chi_1(G)\le{2}$, a contradiction. Thus $(v,u_1)\in{E}$. Similarly
$(v,u_2)\in{E}$.

Now let us assume that $H\cong{P_4}$ and
 $(v,z_2)\in{E(H)}$.
The arguments used to conclude that $v$ and $z$ are both adjacent to
$u_1$ and $u_2$ can now be repeated with reference to the vertices
$z_1$ and $z_2$ since $d_H(z_2)=2$. Thus we conclude, without loss
of generality, that $z_1$ and $z_2$ are both adjacent to say $u_3$
and $u_4$ from $A-\{u_1,u_2\}$. Let $\{u_5\}=A-\{u_1,u_2,u_3,u_4\}$.
Note that $G-u_5\cong{G_2}$.

Now let $H\cong{P_3}\cup{K_1}$. If $z_1$ has at most one neighbour
in $A-\{u_1,u_2\}$ then $\chi_1(G)\le{2}$ since
$$X=A\cup\{z_1\} {\mbox{\rm \  and\ \ }} V-X
{\mbox{\rm \ are both 1-independent}}.$$ Thus $z_1$ and similarly
$z_2$ have at least two neighbours in $A-\{u_1,u_2\}$. Now let
$\{u_3,u_4,u_5\}=A-A_1$. Suppose that $z_1$ and $z_2$ have two
common neighbours in $\{u_3,u_4,u_5\}$, say $u_3$ and $u_4$. Then
clearly $G-u_5\cong{G_1}$.

Now assume that $z_1$ and $z_2$ have exactly one common neighbour.
Specifically, assume that $z_1$ is adjacent to $u_3$ and $u_4$;
$z_2$ is adjacent to $u_3$ and $u_5$. Now
$$X_1=(A-\{u_3\})\cup\{z_1,z_2\}
{\mbox{\rm \  is 1-independent so that \ }} V-X_1 {\mbox{\rm \  is
not\ \ }}
$$
as $\chi_1(G)=3$. This implies that $(v,u_3)\in{E}$. Similarly, by
considering the sets
$$ X_2=\{ u_1,u_2,u_3,u_4,z_2\}
 {\mbox{\rm \  and \ }}
X_3=\{u_1, u_2,u_3,u_5,z_1\}$$ we conclude that $(v,u_5)$ and
$(v,u_4)$ are in $E$. Then $G\cong{G_5}$ given in Figure
~\ref{fig:g5}.

From now onwards we will assume that $H\cong{C_4}$. Thus every
vertex of $H$ has degree $\Delta(H)=2$ in $H$. Moreover we assume
that $z$ has the largest number of neighbours in $A$. Recall that
$(v,z)\notin{E(H)}$. Since $|A_1|\le{2}$, we have
$|N_G(z)\cap{N_G(v)}\cap{A}|\le{2}$.

Firstly if $|N_G(z)\cap{N_G(v)}\cap{A}|= 1$ then the sets
$$X_1=(A-(N_G(z)\cap{N_G(v)}))\cup\{z,v\}
{\mbox{\rm \  and \ }}V-X_1$$ provide a (2,1)-colouring of $G$, a
contradiction to the assumption that $\chi_1(G)=3$.

Next let $|N_G(z)\cap{N_G(v)}\cap{A}|= 0$. If $|A_1| \le {1}$ then
by the choice $z$, $|N_G(v)\cap{A}|\le{1}$. But then the sets
$Y_1=A\cup{\{v,z\}}$ and $V-Y_1=\{u,z_1,z_2\}$ provide a
(2,1)-colouring of $G$, a contradiction. Hence $|A_1|=2$ and let
$A_{1}=\{u_1,u_2\}.$ If $v$ has atmost one neighbour in $A$ then the
sets
$$X_2 = \{v,z,u_2,u_3,u_4,u_5 \}{\mbox{\rm \  and \ }} V-X_2=\{u,u_1,z_1,z_2\}$$
form a (2,1)-colouring of $G$, a contradiction. If $v$ has two
neighbours in $A$, say $u_{3}$ and $u_{4}$, then the sets
$$X_3= \{z_1,z_2,u_1,u_2,u_3,u_4 \}{\mbox{\rm \  and \ }}V-X_3=\{u,u_5,z,v\}$$
provide a (2,1)-colouring of $G$, a contradiction.

Hence $|N_G(z)\cap{N_G(v)}\cap{A}|={2}$. Without any loss of
generality we assume that $N_G(z)\cap{N_G(v)}\cap{A}=\{u_1,u_2\}$.
Similarly we can easily show that
$|N_G(z_1)\cap{N_G(z_2)}\cap{A}|=2$. Without any loss of generality,
let $N_G(z_1)\cap{N_G(z_2)}\cap{A}=\{u_3,u_4\}$. Now let
$\{u_5\}=A-\{u_1,u_2,u_3,u_4\}$. It is easy to see that
$G-u_5\cong{G_3}$.

This completes the proof of the lemma. $\square$

\bigskip


\begin{lemma}
\label{lem:eight} Let $G$ be a triangle-free graph of order 10 with
${\Delta(G)}={4}$ and $3\le\Delta(H)\le{4}$. If $\chi_1(G)={3}$ then
there exists a vertex $u^*$ in $G$ such that $G-u^*\cong{G_1}$ or
$G_2$.
\end{lemma}
\smallskip

\noindent{\it Proof.} We will assume $\chi_1(G)=3$. Let
$A=\{u_1,u_2,u_3,u_4\}$. If $\Delta(H)=4$ then $G$ is a subgraph of
$K_{5,5}$ and $\chi_1(G)\le\chi_0(G)=2$, a contradiction. Hence we
assume $\Delta(H)=3$.

Let $N_H(z)=\{z_1,z_2,z_3\}$ and $v\in{B}$ such that
$(z,v)\not\in{E(H)}$. We provide a proof of this lemma by making and
proving, a sequence of claims.

\noindent{\bf Claim~\ref{lem:eight}.1.} $|N_H(v)|\ge{2}$\\
Suppose that $|N_H(v)|\le{1}$; then we can partition $V$ into two
1-independent sets, $X=A\cup\{z\}$ and $V-X$. Hence
$\chi_1(G)\le{2}$, a contradiction. This establishes
Claim~\ref{lem:eight}.1.

Without any loss of generality, assume that $(v,z_1)$ and $(v,z_2)$
are in $E(H)$. Note that $|N_G(z)\cap{A}|\le{1}$ and
$|N_G(v)\cap{A}|\le{2}$.

\noindent{\bf Claim~\ref{lem:eight}.2.} If $|N_G(z)\cap{A}|=1$ then $G-u_1 \cong{G_2}$.\\
Suppose that $|N_G(z)\cap{A}|=1$ and let $(z,u_1)\in{E}$. If, in
addition, $(v,u_1)\in{E}$ then the sets
$$X=\{u_2,u_3,u_4,z,v\} \ \ {\mbox{\rm \ and\ \ }} V-X$$
partition $V$ into 1-independent sets implying $\chi_1(G)\le{2}$, a
contradiction. Hence $(v,u_1)\not\in{E}$. If $|N_G(v)\cap{A}|\le{1}$
then again $\chi_1(G)\le{2}$ since
$$X_1=A\cup\{v,z\} {\mbox{\rm \  and\ \ }} V-X_1
{\mbox{\rm \  are both 1-independent}},$$ Hence $|N_G(v)\cap{A}|=2$.
Let us assume that $N_G(v)\cap{A}=\{u_2,u_3\}$. The set
$$X_2=\{u_1,u_3,u_4,z,v\} {\mbox{\rm \ is 1 independent, so \ \ }}
V-X_2 {\mbox{\rm \ is not 1-independent}}$$  as $\chi_1(G)=3$. This
implies that $(u_2,z_3)\in{E}$. Similarly we conclude that
$(u_3,z_3)\in{E}$.

\begin{figure}[H]
 \centering
 \setlength{\unitlength}{0.7 cm}
 \begin{picture}(8,4.5)(0.0,1.5)
  \put(0.5,3.5){\circle{0.8}
            \makebox(0,0)[r]{$v$}}
  \put(2,4){%
  \polyline(0,0)(2,0)(2,1.5)(0,1.5)(0,0)}
\put(2.5,5.0){\makebox(0,0)[t]{$z_1$}}
 \put(3.5,5.0){\makebox(0,0)[t]{$z_2$}}
 \put(2.5,4.3){\circle*{0.2}}
 \put(3.5,4.3){\circle*{0.2}}

 \put(5,4){%
  \polyline(0,0)(2,0)(2,1.5)(0,1.5)(0,0)}
\put(5.5,5.0){\makebox(0,0)[t]{$z$}}
 \put(6.5,5.0){\makebox(0,0)[t]{$u_4$}}
 \put(5.5,4.3){\circle*{0.2}}
 \put(6.5,4.3){\circle*{0.2}}
 \put(2,1.5){%
  \polyline(0,0)(2,0)(2,1.5)(0,1.5)(0,0)}
\put(2.5,1.8){\makebox(0,0)[b]{$u_2$}}
 \put(3.5,1.8){\makebox(0,0)[b]{$u_3$}}
 \put(2.5,2.5){\circle*{0.2}}
 \put(3.5,2.5){\circle*{0.2}}
 \put(5,1.5){%
  \polyline(0,0)(2,0)(2,1.5)(0,1.5)(0,0)}
\put(5.5,1.8){\makebox(0,0)[b]{$u$}}
 \put(6.5,1.8){\makebox(0,0)[b]{$z_3$}}
 \put(5.5,2.5){\circle*{0.2}}
 \put(6.5,2.5){\circle*{0.2}}

\polyline(5.5,4.3)(6.5,2.5)(6.5,4.3)(5.5,2.5)
 \polyline(0.85,3.73)(2,4.75)
  \polyline(0.9,3.58)(2,4.60)
  \polyline(0.9,3.4)(2,2.6)
  \polyline(0.85,3.25)(2,2.45)

 \polyline(4,4.75)(5.0,4.75)
  \polyline(4,4.6)(5.0,4.6)
  \polyline(4,2.6)(5.0,2.6)
  \polyline(4,2.45)(5.0,2.45)
 \end{picture}
\caption{\label{fig:Gu1} $G-u_1\cong{G_2}$}
\end{figure}
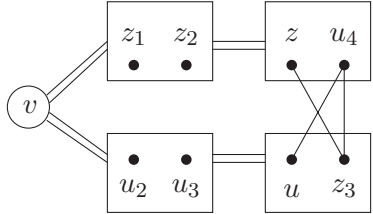

Since the sets
$$Y_1=\{u,z,v,z_3\} {\mbox{\rm \  and\ \ }} Y_2=\{u_1,u_2,u_3,z_1,z_2\}{\mbox{\rm \  are  1-independent,}}$$
  $$V-Y_1=A\cup\{z_1,z_2\}{\mbox{\rm \  and\ \ }}V-Y_2=\{z,z_3,u,u_4,v\}{\mbox{\rm \ are not 1-independent}}$$
as $\chi_1(G)={3}$. Hence $(u_4,z_1)$, $(u_4,z_2)$ and $(u_4,z_3)$
are all in E. Now $G-u_1$ is isomorphic to $G_2$ given in Figure
~\ref{fig:Gu1}.
\smallskip

\noindent
 This establishes Claim~\ref{lem:eight}.2. Henceforth we
will assume that $|N_G(z)\cap{A}|=0$.

\noindent{\bf Claim~\ref{lem:eight}.3.}  $|N_G(v)\cap{A}|=2$ and $(v,z_3)\not\in{E(H)}$.\\
Otherwise, that is, if $|N_G(v)\cap{A}|\le{1}$, then
$X=A\cup\{z,v\}$ and $V-X$ provide a partition of $V$ into
1-independent sets, implying $\chi_1(G)\le{2}$. Hence
$|N_G(v)\cap{A}|=2$. Since $d_G(v)\le{4}$ we now have
$(v,z_3)\not\in{E}$. This establishes Claim~\ref{lem:eight}.3.

Without any loss of generality, we now assume that
$N_G(v)\cap{A}=\{u_1,u_2\}$. Clearly there are no edges between
$\{z_1,z_2\}$ and $\{u_1,u_2\}$.

\noindent{\bf Claim~\ref{lem:eight}.4.}  For $i=1$ and 2, $(u_i,z_3)\in{E}$.\\
Now note that the set $ X_1= \{u_2,u_3,u_4,z,v\}$  is 1-independent
 while  $V-X_1 = \{u,u_1,z_1,z_2,z_3\}$ is not
as $\chi_1(G)=3$. This implies  $(u_1,z_3)\in{E}$. Similarly
$(u_2,z_3)\in{E}$. This establishes Claim~\ref{lem:eight}.4.

\begin{figure}[H]
 \centering
  \setlength{\unitlength}{0.7 cm}
 \begin{picture}(4.5,4.5)(0.0,0.5)
  \put(0.25,3.5){\circle{0.8}
            \makebox(0,0)[r]{$u$}}
  \put(1.5,4){%
  \polyline(0,0)(1.5,0)(1.5,1.5)(0,1.5)(0,0)}
\put(1.75,5.0){\makebox(0,0)[t]{$u_1$}}
 \put(2.75,5.0){\makebox(0,0)[t]{$u_2$}}
 \put(1.75,4.3){\circle*{0.2}}
 \put(2.75,4.3){\circle*{0.2}}
 \put(4,4){%
  \polyline(0,0)(1.5,0)(1.5,1.5)(0,1.5)(0,0)}
\put(4.25,5.0){\makebox(0,0)[t]{$v$}}
 \put(5.25,5.0){\makebox(0,0)[t]{$z_3$}}
 \put(4.25,4.3){\circle*{0.2}}
 \put(5.25,4.3){\circle*{0.2}}
 \put(1.5,1.5){%
  \polyline(0,0)(1.5,0)(1.5,1.5)(0,1.5)(0,0)}
\put(1.75,1.8){\makebox(0,0)[b]{$u_3$}}
 \put(2.75,1.8){\makebox(0,0)[b]{$u_4$}}
 \put(1.75,2.5){\circle*{0.2}}
 \put(2.75,2.5){\circle*{0.2}}
 \put(4,1.5){%
  \polyline(0,0)(1.5,0)(1.5,1.5)(0,1.5)(0,0)}
\put(4.25,1.8){\makebox(0,0)[b]{$z_1$}}
 \put(5.25,1.8){\makebox(0,0)[b]{$z_2$}}
 \put(4.25,2.5){\circle*{0.2}}
 \put(5.25,2.5){\circle*{0.2}}

\polyline(4.25,2.5)(4.25,4.3)(5.25,2.5)
 \polyline(0.55,3.73)(1.5,4.75)
  \polyline(0.65,3.65)(1.5,4.60)
  \polyline(0.65,3.35)(1.5,2.6)
  \polyline(0.55,3.23)(1.5,2.45)
 \polyline(3,4.75)(4,4.75)
  \polyline(3,4.6)(4,4.6)
  \polyline(3,2.6)(4,2.6)
  \polyline(3,2.45)(4,2.45)
 \put(2.5,1.0){\makebox(0,0)[l]{$(a)$ $G-z\cong{G_1}$}}
 \end{picture}
\hspace{1.5cm}
 \begin{picture}(4.5,4.5)(0.0,0.5)
  \put(0.25,3.5){\circle{0.8}
            \makebox(0,0)[r]{$v$}}
  \put(1.5,4){%
  \polyline(0,0)(1.5,0)(1.5,1.5)(0,1.5)(0,0)}
\put(1.75,5.0){\makebox(0,0)[t]{$u_1$}}
 \put(2.75,5.0){\makebox(0,0)[t]{$u_2$}}
 \put(1.75,4.3){\circle*{0.2}}
 \put(2.75,4.3){\circle*{0.2}}

 \put(4,4){%
  \polyline(0,0)(1.5,0)(1.5,1.5)(0,1.5)(0,0)}
\put(4.25,5.0){\makebox(0,0)[t]{$u$}}
 \put(5.25,5.0){\makebox(0,0)[t]{$z_3$}}
 \put(4.25,4.3){\circle*{0.2}}
 \put(5.25,4.3){\circle*{0.2}}
 \put(1.5,1.5){%
  \polyline(0,0)(1.5,0)(1.5,1.5)(0,1.5)(0,0)}
\put(1.75,1.8){\makebox(0,0)[b]{$z_1$}}
 \put(2.75,1.8){\makebox(0,0)[b]{$z_2$}}
 \put(1.75,2.5){\circle*{0.2}}
 \put(2.75,2.5){\circle*{0.2}}
 \put(4,1.5){%
  \polyline(0,0)(1.5,0)(1.5,1.5)(0,1.5)(0,0)}
\put(4.25,1.8){\makebox(0,0)[b]{$z$}}
 \put(5.25,1.8){\makebox(0,0)[b]{$u_4$}}
 \put(4.25,2.5){\circle*{0.2}}
 \put(5.25,2.5){\circle*{0.2}}

\polyline(4.25,4.3)(5.25,2.5)(5.25,4.3)(4.25,2.5)
 \polyline(0.55,3.73)(1.5,4.75)
  \polyline(0.65,3.65)(1.5,4.60)
  \polyline(0.65,3.35)(1.5,2.6)
  \polyline(0.55,3.23)(1.5,2.45)

 \polyline(3,4.75)(4,4.75)
  \polyline(3,4.6)(4,4.6)
  \polyline(3,2.6)(4,2.6)
  \polyline(3,2.45)(4,2.45)
\put(2.5,1.0){\makebox(0,0)[l]{$(b)$ $G-u_3\cong{G_2}$}}
 \end{picture}
\caption{\label{fig:Gz} $Graph$  $G-u^*$}
\end{figure}
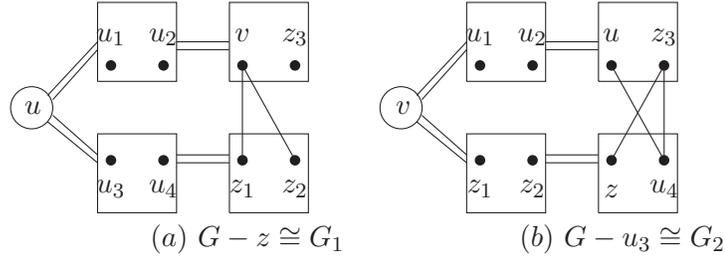

Since $z_3$ is adjacent to $u_1$, $u_2$ and $z$ and $d_G(z_3)\le{4}$
we can assume, without any loss of generality, that
$(z_3,u_3)\not\in{E}$. The set $ X_1= \{u,u_3,v,z,z_3\}$  is
1-independent while  $V-X_1 = \{u_1,u_2,u_4,z_1,z_2\}$ cannot be as
$\chi_1(G)=3$. This implies that $(u_4,z_1)$ and $(u_4,z_2)$ are
both in $E$. Now if $(z_3,u_4)\not\in{E}$, we can similarly conclude
that $(u_3,z_i)\in{E}$ for $i=1$ and 2. In this case we can easily
verify that $G-z\cong{G_1}$ (see Figure ~\ref{fig:Gz}$(a)$). On the
other hand, that is if $(z_3,u_4)\in{E}$, we can check that
$G-u_3\cong{G_2}$ (see Figure ~\ref{fig:Gz}$(b)$).
\smallskip

\noindent

This proves the lemma. $\square$

\medskip

Suppose that $G$ is a triangle-free graph of order 10 with
$\Delta(G)=4$ and $\chi_1(G)=3$. As a consequence of
Lemmas~\ref{lem:five}(i) and \ref{lem:eight} we can assume that
$\Delta(H)=2$. It is easy to see that $H$ is isomorphic to one of
the graphs (i) $P_3\cup{2K_1}$
 (ii) $P_3\cup{K_2}$
(iii) $P_4\cup{K_1}$ (iv) $P_5$ (v) $C_5$ and (vi) $C_4\cup{K_1}$.

\medskip

\begin{lemma}
\label{lem:nine} Let $G$ be a  triangle-free graph of order 10 with
${\Delta(G)}={4}$ and $\Delta(H)=2$. Furthermore, let $H$ be
isomorphic to $P_3\cup{2K_1}$ or $P_3\cup{K_2}$. If $\chi_1(G)={3}$
then there exists a vertex $u^*$ in $G$ such that $G-u^*\cong{G_1}$
or $G_2$ or $G_3$.
\end{lemma}
\smallskip

\noindent{\it Proof.} Assume that $\chi_1(G)=3$. Let $z\in{B}$ with
$d_H(z)=2$ and $N_H(z)=\{z_1,z_2\}$. For $x$ in $\{ z, z_1,  z_2 \}$
we have $|N_G(x)\cap{A}|\ge{2}$, otherwise $X_1=A\cup\{x\}$ and
$V-X_1$ provide a $(2,1)$-colouring of $G$, a contradiction. Since
$d_H(z)=2$, we have $|N_G(z)\cap{A}|={2}$. Since $G$ is $K_3$-free,
this implies $|N_G(z_i)\cap{A}|={2}$ for $i=1$ and 2. Without any
loss of generality we can write $N_G(z)\cap{A}=\{u_1,u_2\}$ and
$N_G(z_i)\cap{A}=\{u_3,u_4\}$ for $i=1$ and 2.

Let $\{z_3,z_4\}=V(H)-\{z,z_1,z_2\}$. If $(z_3,u_1)$ and $(z_3,u_2)$
are in $E$ then $G-z_4\cong{G_1}$ or $G_2$ or $G_3$ according as the
number of edges between $\{z_3\}$ and $\{u_3,u_4\}$ is 0 or 1 or 2.
Hence we will assume, without loss of generality, that
$(z_3,u_2)\not\in{E}$. Suppose $(z_4,u_2)\not\in{E}$ then
$$X_1=\{u,u_2,z_1,z_2,z_3,z_4\}
{\mbox{\rm \ and\ \ }}  V-X_1
$$
form a (2,1)-colouring of $G$, a contradiction. Hence
$(z_4,u_2)\in{E}$. If $(z_4,u_1)\in{E}$ then $G-z_3\cong{G_1}$ or
$G_2$ or $G_3$. Hence we assume that $(z_4,u_1)\not\in{E}$. Now
since $d_G(u_3)\le{4}$, we can assume that $(u_3,z_3)\not\in{E}$,
from which it follows that the sets
$$X_1=\{u_2,u_3,u_4,z,z_3\}
{\mbox{\rm \ and\ \ }}  V-X_1
$$
form a (2,1)-colouring of $G$, a contradiction.

This proves the lemma. $\square$

\medskip

\begin{lemma}
\label{lem:ten} Let $G$ be a  triangle-free  graph of order 10 with
${\Delta(G)}={4}$ and $\Delta(H)=2$. Furthermore, let $H$ be
isomorphic to $P_4\cup{K_1}$. If $\chi_1(G)={3}$ then there exists a
vertex $u^*$ in $G$ such that $G-u^*\cong{G_i}$, for some $i$, $1
\le {i} \le {3}$.
\end{lemma}
\smallskip

\noindent{\it Proof.} Let us suppose that $\chi_1(G)=3$.  Let $z$
and $z_1$ be vertices in $B$ with $d_H(z)=d_H(z_1)=2$. Note that
$(z,z_1)\in{E(H)}$. Let $z_2$ ($z_3$) be the other neighbour of $z$
($z_1$) . Finally, let $\{z_4\}=V(H)-\{z,z_1,z_2,z_3\}$.

\noindent{\bf Claim~\ref{lem:ten}.1.}
For $x=z$ and $z_1$, $|N_G(x)\cap{A}|=2$.\\
This claim can be proved using arguments similar to the ones used in
Lemma~\ref{lem:nine}.

Now, without any loss of generality, let $N_G(z)\cap{A}=\{u_1,u_2\}$
and $N_G(z_1)\cap{A}=\{u_3,u_4\}$. Since $\chi_1(G)=3$ and
$V-A-\{z_2,z_3\}$ is 1-independent it follows that
$A\cup\{z_2,z_3\}$ is not 1-independent. Note that $z_2$ and $z_3$
do not have a common neighbour in $A$. Thus we conclude that either
$(z_2,u_i)\in{E}$ for $i=3$ and 4 or $(z_3,u_i)\in{E}$ for $i=1$ and
2. Suppose, without loss of generality, $(z_2,u_i)\in{E}$ for $i=3$
and 4.

If $z_3$ is adjacent to both $u_1$ and $u_2$, then it is easy to
verify that $G-z_4\cong{G_2}$.

Hence $(z_3,u_i)\not\in{ E}$ for $i=1$ or 2. Without any loss of
generality assume that $(z_3,u_1)\not\in{E}$. Now
$$X_1=\{u_2,u_3,u_4,z\} {\mbox{\rm \ and \ \ }} X_2=\{u_1,u_3,u_4,z,z_3\}
{\mbox{\rm \ are 1-independent. \ \ }}$$ Since $\chi_1(G)=3$, the
sets $V-X_1=\{u,u_1,z_1,z_2,z_3,z_4\}$  and
$V-X_2=\{u,u_2,z_1,z_2,z_4\}$ are not 1-independent. This in turn
implies that $(u_i,z_4)\in{E}$ for $i=$ 1 and 2.  Now it is easy to
verify that  $G-z_3\cong{G_1{\mbox{\rm \ or }}G_2 {\mbox{\rm \ or }}
G_3}$.

Hence it follows that there exists a $u^*$ such that
$G-u^*\cong{G_1{\mbox{\rm \ or\ \ }} G_2 {\mbox{\rm \ or\ \ }}G_3}$.

This proves the lemma. $\square$
\medskip

\begin{lemma}
\label{lem:eleven} Let $G$ be a triangle-free graph of order 10 with
${\Delta(G)}={4}$ and $\Delta(H)=2$. Furthermore, let $H$ be
isomorphic to $P_5$. If $\chi_1(G)={3}$ then there exists a vertex
$u^*$ such that $G-u^*\cong{G_i}$, for some $i$, $1 \le {i} \le
{3}$.
\end{lemma}
\smallskip

\noindent{\it Proof.}
 We assume that $\chi_1(G)=3$. Let $z$ be the central
vertex of $H$. Since $\Delta(G)=4$, $|N_G(z)\cap{A}|\le{2}$. If
$|N_G(z)\cap{A}|\le{1}$ then $X=A\cup\{z\}$ and $V-X$ form a
partition of $V$ into 1-independent sets implying $\chi_1(G)\le{2}$.
Thus $|N_G(z)\cap{A}|={2}$ and let $N_G(z)\cap{A}=\{u_1,u_2\}$. Also
let $N_H(z)=\{z_1,z_2\}$. Furthermore, let $z_3$ and $z_4$ be the
neighbours of $z_1$ and $z_2$ respectively.

Since $\chi_1(G)=3$ and
 $X=\{u,z,z_3,z_4\}$ is 0-independent, the set  $V-X=A\cup\{z_1,z_2\}$
 is not 1-independent.

Since $\{u_1,u_2,z_1,z_2\}$ is totally disconnected, it follows that
$\Delta(L)=2$ where $L=G[\{u_3,u_4,z_1,z_2\}]$. Suppose that
$d_L(u_3)=\Delta(L)=2$. This means that $(u_3,z_i)\in{E}$ for $i=1$
and 2. Since $G$ is triangle-free $(u_3,z_i)\not\in{E}$  for $i=3$
and 4.

Now note that $d_G(z)=\Delta(G)=4$. Let
$$F=G[V-N_G[z]]
=G[\{u,u_3,u_4,z_3,z_4\}].$$
Clearly either\\
(i) $\Delta(F)=d_F(u_4)=3$ or\\
(ii) $\Delta(F)=2$ and $F\cong{P_4}\cup{K_1}$ or $P_3\cup{2K_1}$.\\
Hence Lemma~\ref{lem:eleven} is established using
Lemmas~\ref{lem:eight} to~\ref{lem:ten} in the case
$d_L(u_3)=\Delta(L)=2$.
Similarly, the lemma is established when $d_L(u_4)=\Delta(L)=2$: in
other words when $(u_4,z_i)\in{E}$ for $i=1$ and 2.

Now let us assume that $d_L(z_1)=\Delta(L)=2$, that is
$(z_1,u_i)\in{E}$ for $i=3$ and 4. Therefore $(z_3,u_i)\not\in{E}$
for $i=3$ and $4$. Now note that $d_G(z)=\Delta(G)=4$.

Note that
$F=G[V-N_G[z]]=G[\{u,u_3,u_4,z_3,z_4\}]\cong{P_3\cup{2K_1}}$ or
$P_4\cup{K_1}$ or $C_4\cup{K_1}$ according as $z_4$ is adjacent to 0
or 1 or 2 vertices from $\{u_3,u_4\}$.

If $F\cong{P_3\cup{2K_1}}$  or $P_4\cup{K_1}$ then
Lemma~\ref{lem:eleven} is established using Lemmas~\ref{lem:nine}
and \ref{lem:ten}.

Hence we assume that $F\cong{C_4\cup{K_1}}$. This implies that
$(z_4,u_i)$ is in E for $i=$ 3 and 4. Since $\chi_1(G)=3$ and the
set
$$X_1=\{u,z,z_1,z_4\}
{\mbox{\rm \ is 1-independent, the set  \ \ }} V-X_1 {\mbox{\rm \ is
not 1-independent}}. $$

Thus $(z_3,u_i)\in{E}$ for $i=1,~2$. Now it is easy to verify that
$G-z_2\cong{G_1}$ or $G_2$ according as the number of edges between
$\{z_4\}$ and $\{u_1, u_2\}$ is 0 or 1. This establishes the lemma
when $d_L(z_1)=\Delta(L)=2$.

Since the vertices $z_1$ and $z_2$ are similar, the lemma is
established when $d_L(z_2)=\Delta(L)=2$ in a similar manner.

This completes the proof of Lemma~\ref{lem:eleven}. $\square$

\medskip

\begin{lemma}
\label{lem:twelve} Let $G$ be a triangle-free graph of order 10 with
${\Delta(G)}={4}$ and $\Delta(H)=2$. Furthermore, let $H$ be
isomorphic to $C_5$. If $\chi_1(G)={3}$ then there exists a vertex
$u^*$ in $G$ such that $G-u^*\cong{G_i}$ for some $i$,
$1\le{i}\le{3}$.
\end{lemma}
\smallskip

\noindent{\it Proof.} Let $V(H)=\{z_1,z_2,z_3,z_4,z_5\}$. Assume
that $(z_i,z_{i+1})\in{E(H)}$ for $i=1$, 2, 3, 4 and
$(z_5,z_1)\in{E(H)}$. Assume that $\chi_1(G)=3$. The set $
X_1=\{u,z_2,z_4,z_5\}$ is 1-independent and so
$V-X_1=A\cup\{z_1,z_3\}$  is not 1-independent.

This implies that $\Delta(L)=2$ where $L=G[A\cup\{z_1,z_3\}]$. Now,
either $d_L(u_i)=2$ for some $i$, $1\le{i}\le{4}$ or
$d_L(z_i)=2$ for $i=1$ or 3.\\
{\bf Case i.} 
$d_L(u_i)=2$ for some $i$, say $i=1$.

Hence $(u_1,z_i)\in{E}$ for $i=1$ and 3. Since $G$ is triangle-free,
$(u_1,z_i)\not\in{E}$ for $i=2$, 4, 5. Since $\chi_1(G)=3$ and the
set $Y_1=\{u,u_1,z_2,z_4,z_5\}$ is $1-$independent,  the set
$V-Y_1=\{u_2,u_3,u_4,z_1,z_3\}$ is not 1-independent. This in turn
implies that, for some $i\in\{2,3,4\}$, $(u_i,z_j)\in{E}$ for $j=1$
and 3. Without any loss of generality we assume that
$(u_2,z_j)\in{E}$ for $j=1$ and 3. Now note that
$(u_2,z_j)\not\in{E}$ for $j=2$, 4 and 5. Observe that
$d_G(z_1)=\Delta(G)=4$. Let
$F=G[V-N_G[z_1]]=G[\{u,u_3,u_4,z_3,z_4\}]$.
Clearly either\\
(i) $\Delta(F)=3$, or\\
(ii) $F\cong{P_3}\cup{K_2}$ or $P_5$.\\
Thus Lemma~\ref{lem:twelve} is established using
Lemmas~\ref{lem:eight}, \ref{lem:nine} and \ref{lem:eleven}, in Case
i.

\noindent
{\bf Case ii.} 
$d_L(z_i)=2$ for $i=1$ or 3.\\
Let us assume that $(z_1,u_i)\in{E}$ for $i=1$ and 2. Note that
$d_G(z_1)=4$ and consider the subgraph
$G[V-N_G[z_1]]=G[\{u,u_3,u_4,z_3,z_4\}]=F$, say. Since $G$ is
triangle-free, the vertex $u_3$ (also $u_4$) is adjacent to at most
one of $z_3$ and $z_4$. If $u_3$ (or $u_4$) is adjacent to neither
$z_3$ nor $z_4$ then $F\cong{P_3}\cup{K_2}$ or $P_5$. Thus the lemma
is established using Lemmas~\ref{lem:nine} and~\ref{lem:eleven}.
Suppose that both $u_3$ and $u_4$ are adjacent to the same vertex,
say $z_3$, then $\Delta(F)=3$ and the lemma is established using
Lemma~\ref{lem:eight}. Hence without any loss of generality assume
that $(u_3,z_3)$ and $(u_4,z_4)$ are in $E$. Hence $(u_3,z_2)$ and
$(u_4,z_5)$ are not in $E$. Now, it is easy to check that
$$Y_1=\{u_1,u_2,u_3,z_2,z_4\}
{\mbox{\rm \ and\ }} V-Y_1=\{u,u_4,z_1,z_3,z_5\} $$ provide a
(2,1)-colouring of $G$, a contradiction. This proves
Lemma~\ref{lem:twelve}. $\square$

\medskip

\begin{lemma}
\label{lem:thirteen} Let $G$ be a triangle-free graph of order 10
with ${\Delta(G)}={4}$ and $\Delta(H)=2$. Furthermore, let $H$ be
isomorphic to $C_4\cup{K_1}$. If $\chi_1(G)={3}$ then there exists a
vertex $u^*$ such that $G-u^*\cong{G_i}$ for some $i$,
$1\le{i}\le{3}$.
\end{lemma}
\smallskip

\noindent{\it Proof.} Let us assume that $\chi_1(G)=3$. Recall that
$u\in{V}$ with $d_G(u)=\Delta(G)=4$, $N_G(u)=A=\{u_1,u_2,u_3,u_4\}$,
$B=\{z_1,z_2,z_3,z_4,z_5\}$ and $H=G[B]=C_4\cup{K_1}$. Assume that
$(z_i,z_{i+1})\in{E(H)}$ for $i=1$, 2, 3 and $(z_4,z_1)\in{E(H)}$.
Hence $z_5$ has degree 0 in $H$.

The sets
$$Y_1=\{u,z_2,z_4,z_5\}\ \
{\mbox{\rm and\ }} Y_2=\{u,z_1,z_3,z_5\} \ \ {\mbox{\rm are
1-independent.\ }}$$ Since $\chi_1(G)=3$ the sets
$$V-Y_1=\{z_1,z_3\}\cup{A}\ \
{\mbox{\rm and\ }} V-Y_2=\{z_2,z_4\}\cup{A} \ \ {\mbox{\rm are not
1-independent.\ }}$$ Hence $F_1=G[V-Y_1]$ and $F_2=G[V-Y_2]$ both
have maximum degree 2.

\noindent{\bf Case i.} 
The subgraph $F_i$, $i=1,\ 2$, attains its maximum degree at a $z_j$
for some $j$ in $\{1,2,3,4\}$. We assume without loss of generality
that
$$d_{F_1}(z_1)=2,\
N_{F_1}(z_1)=\{u_1,u_2\},\ \ d_{F_2}(z_2)=2,\
N_{F_2}(z_2)=\{u_3,u_4\}.$$

\begin{figure}[H]
 \centering
  \setlength{\unitlength}{0.7 cm}
 \begin{picture}(4.5,4.5)(0.0,0.5)
  \put(0.25,3.5){\circle{0.8}
            \makebox(0,0)[r]{$u_1$}}
  \put(1.3,4){%
  \polyline(0,0)(1.5,0)(1.5,1.5)(0,1.5)(0,0)}
\put(1.55,5.0){\makebox(0,0)[t]{$z_1$}}
 \put(2.55,5.0){\makebox(0,0)[t]{$z_3$}}
 \put(1.55,4.3){\circle*{0.2}}
 \put(2.55,4.3){\circle*{0.2}}
 \put(3.8,4){%
  \polyline(0,0)(1.5,0)(1.5,1.5)(0,1.5)(0,0)}
\put(4.05,5.0){\makebox(0,0)[t]{$z_2$}}
 \put(5.05,5.0){\makebox(0,0)[t]{$z_4$}}
 \put(4.05,4.3){\circle*{0.2}}
 \put(5.05,4.3){\circle*{0.2}}
 \put(1.3,1.5){%
  \polyline(0,0)(1.5,0)(1.5,1.5)(0,1.5)(0,0)}
\put(1.55,1.8){\makebox(0,0)[b]{$u$}}
 \put(2.55,1.8){\makebox(0,0)[b]{$z_5$}}
 \put(1.55,2.5){\circle*{0.2}}
 \put(2.55,2.5){\circle*{0.2}}
 \put(3.8,1.5){%
  \polyline(0,0)(1.5,0)(1.5,1.5)(0,1.5)(0,0)}
\put(4.05,1.8){\makebox(0,0)[b]{$u_3$}}
 \put(5.05,1.8){\makebox(0,0)[b]{$u_4$}}
 \put(4.05,2.5){\circle*{0.2}}
 \put(5.05,2.5){\circle*{0.2}}

\polyline(4.05,2.5)(4.05,4.3)(5.05,2.5)
 \Dline(4.05,2.5)(5.05,4.3){0.2}
 \Dline(5.05,4.3)(5.05,2.5){0.2}
 \polyline(0.55,3.73)(1.3,4.75)
  \polyline(0.65,3.65)(1.3,4.60)
  \polyline(0.65,3.35)(1.3,2.6)
  \polyline(0.55,3.23)(1.3,2.45)
 \polyline(2.8,4.75)(3.8,4.75)
  \polyline(2.8,4.6)(3.8,4.6)
  \polyline(2.8,2.6)(3.8,2.6)
  \polyline(2.8,2.45)(3.8,2.45)
 \put(0.8,1.0){\makebox(0,0)[l]{$(a)$ $G-u_2\cong{G_i}$, $1\le i \le 3$}}
 \end{picture}
\hspace{2.0cm}
 \begin{picture}(4.5,4.5)(0.0,0.5)
  \put(0.0,2.5){%
  \polyline(0,0)(0.8,0)(0.8,2)(0,2)(0,0)}
  \put(0.35,4.2){\makebox(0,0)[t]{$u_2$}}
  \put(0.35,3.6){\circle*{0.2}}
  \put(0.35,2.8){\makebox(0,0)[b]{$(u_3)$}}

  \put(1.3,4){%
  \polyline(0,0)(1.5,0)(1.5,1.5)(0,1.5)(0,0)}
\put(1.55,5.0){\makebox(0,0)[t]{$u$}}
 \put(2.55,5.0){\makebox(0,0)[t]{$z_5$}}
 \put(1.55,4.3){\circle*{0.2}}
 \put(2.55,4.3){\circle*{0.2}}

 \put(3.8,4){%
  \polyline(0,0)(1.5,0)(1.5,1.5)(0,1.5)(0,0)}
\put(4.05,5.0){\makebox(0,0)[t]{$u_1$}}
 \put(5.05,5.0){\makebox(0,0)[t]{$u_4$}}
 \put(4.05,4.3){\circle*{0.2}}
 \put(5.05,4.3){\circle*{0.2}}
 \put(1.3,1.5){%
  \polyline(0,0)(1.5,0)(1.5,1.5)(0,1.5)(0,0)}
\put(1.55,1.8){\makebox(0,0)[b]{$z_2$}}
 \put(2.55,1.8){\makebox(0,0)[b]{$z_4$}}
 \put(1.55,2.5){\circle*{0.2}}
 \put(2.55,2.5){\circle*{0.2}}
 \put(3.8,1.5){%
  \polyline(0,0)(1.5,0)(1.5,1.5)(0,1.5)(0,0)}
\put(4.05,1.8){\makebox(0,0)[b]{$z_1$}}
 \put(5.05,1.8){\makebox(0,0)[b]{$z_3$}}
 \put(4.05,2.5){\circle*{0.2}}
 \put(5.05,2.5){\circle*{0.2}}

\polyline(4.05,2.5)(4.05,4.3)(5.05,2.5)
 \polyline(0.8,3.75)(1.3,4.75)
  \polyline(0.8,3.60)(1.3,4.60)
  \polyline(0.8,3.6)(1.3,2.6)
  \polyline(0.8,3.45)(1.3,2.45)

 \polyline(2.8,4.75)(3.8,4.75)
  \polyline(2.8,4.6)(3.8,4.6)
  \polyline(2.8,2.6)(3.8,2.6)
  \polyline(2.8,2.45)(3.8,2.45)
\put(0.0,1.0){\makebox(0,0)[l]{$(b)$ $G-u_3\cong{G_1}$ (resp.
$G-u_2\cong{G_1}$)}}
 \end{picture}
\caption{\label{fig:thirteenGu2} $Graph$  $G-u^*$}
\end{figure}
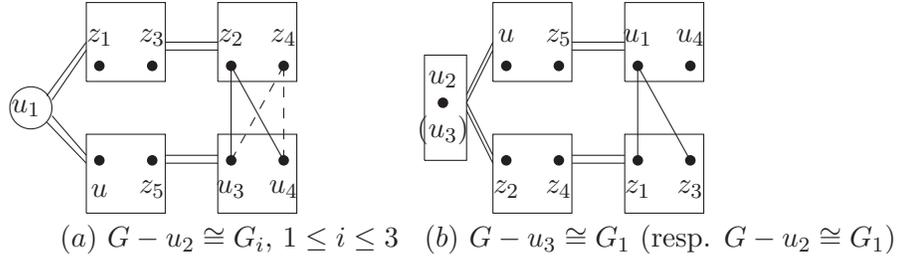

Note that $d_G(z_i)=4$ for $i=1$ and 2. Now we can assume that the
subgraphs $L_1=G[V-N_G[z_1]] =G[\{u,u_3,u_4,z_3,z_5\}]$ and
$L_2=G[V-N_G[z_2]] =G[\{u,u_1,u_2,z_4,z_5\}]$ are both isomorphic to
$C_4\cup{K_1}$. For otherwise by Lemmas~\ref{lem:eight} to
\ref{lem:twelve} there exists a vertex $u^*$ in $G$ such that
$G-u^*\cong{G_i}$ for some $i$, $1\le{i}\le{3}$. Thus
$(z_5,u_i)\in{E}$ for $i=1$, 2, 3 and 4. Now the set
$$X_1=\{z_1,z_2,z_5,u\} \ \
{\mbox{\rm is 1-independent and so \ }} V-X_1=A\cup{\{z_3,z_4\}} \ \
{\mbox{\rm is not}}.$$ Hence we can assume, without loss of
generality, that $(z_3,u_1)\in{E}$. It is easy to verify that the
graph $G-u_2\cong{G_1}$ or $G_2$ or $G_3$ according as the number of
edges between $z_4$ and $\{u_3,u_4\}$ is 0 or 1 or 2. The graph
$G-u_2$ is illustrated in Figure~\ref{fig:thirteenGu2}$(a)$. The
dotted lines indicate that the edges may or may not be in $G$. This
completes the proof of Lemma~\ref{lem:thirteen} in
Case i. 


\noindent

{\bf Case ii.} 
The subgraph $F_1$ attains its maximum degree at a $u_j$ for some
$j$ in $\{1,2,3,4\}$ and $F_2$ attains its maximum degree at a $z_j$
for some $j$ in $\{2,4\}$. Furthermore $d_{F_1}(z_i)\le{1}$ for
$i=1$ and 3.

We assume without loss of generality that $(u_1,z_i)\in{E(F_1)}$ for
$i=1$ and 3; $(u_j,z_2)\in{E(F_2)}$ for $j=2$ and 3. Note that
$N_{F_1}(z_1)=N_{F_1}(z_3)=\{u_1\}$.

Since $d_G(z_2)=4$, the subgraph
$$M_1=G[V-N_G[z_2]]
=G[\{u,u_1,u_4,z_4,z_5\}] $$ can be assumed to be isomorphic to
$C_4\cup{K_1}$. For otherwise, the lemma is established using
Lemmas~\ref{lem:eight} to \ref{lem:twelve}. Hence $(z_5,u_i)\in{E}$
for $i=1$ and 4 and $(z_4,u_4)\not\in{E}$. Thus $d_G(u_1)=4$. Again
$$M_2=G[V-N_G[u_1]]
=G[\{u_2,u_3,u_4,z_2,z_4\}] $$ is assumed to be isomorphic to
$C_4\cup{K_1}$,  by Lemmas~\ref{lem:eight} to \ref{lem:twelve}.
Hence $(z_4,u_2)$ and $(z_4,u_3)$ are in $E$. The set
$X_1=\{u,u_1,z_2,z_4\}$ is 1-independent and so
$V-X_1=\{u_2,u_3,u_4,z_1,z_3,z_5\}$ is not as $\chi_1(G)=3$. This
implies that $z_5$ is adjacent to at least one of $\{u_2,u_3\}$. If
$z_5$ is adjacent to $u_2$ (resp. $u_3$) then it is easy to check
that $G-u_3\cong{G_1}$ (resp. $G-u_2\cong{G_1}$). The graph $G-u_3$
(resp. $G-u_2$) is illustrated in Figure~\ref{fig:thirteenGu2}$(b)$.
This completes the proof of Lemma~\ref{lem:thirteen} in
Case ii. 

\noindent{\bf Case iii.} 
Each subgraph $F_i$, $i=1,\ 2$, attains its maximum degree at a
$u_j$ for some $j$ in $\{1,2,3,4\}$. Furthermore, every $z_j$ has
degree at most 1 in the corresponding $F_i$. We assume without loss
of generality that
$$d_{F_1}(u_1)=2,\
N_{F_1}(u_1)=\{z_1,z_3\},\ \ d_{F_2}(u_2)=2,\
N_{F_2}(u_2)=\{z_2,z_4\}.$$ Note that there are no other edges
between $A$ and $\{z_1,z_2,z_3,z_4\}$. The set
$X_1=\{u_2,u_3,u_4,z_1,z_3\}$ is 1-independent and so
$V-X_1=\{u,u_1,z_2,z_4,z_5\}$ is not as $\chi_1(G)=3$. Hence
$(z_5,u_1)\in{E}$. Now note that $d_G(u_1)=4$. But
$$N_1=G[V-N_G[u_1]]
=G[\{u_2,u_3,u_4,z_2,z_4\}] $$ is isomorphic to $P_3\cup{2K_1}$.
Hence by Lemma~\ref{lem:nine} there exists a vertex $u^*$ such that
$G-u^*\cong{G_i}$ for some $i$, $1\le{i}\le{3}$.

This completes the proof of the lemma. $\square$

\medskip

Combining Lemmas~\ref{lem:five} to \ref{lem:thirteen} we have the
following result.

\begin{theorem}
\label{thm:four} Let $G$ be a triangle-free  graph of order 10 with
$\chi_1(G)=3$. Then either $G\cong{G_5}$ given in
Figure~\ref{fig:g5} or there exists a vertex $u^*$ such that
$G-u^*\cong{G_i}$ for some $i$, $1\le{i}\le{4}$.
\end{theorem}

We observe that there are exactly four triangle-free graphs of order
9, namely $G_i$, $1\le{i}\le{4}$ which are $(3,1)$-critical.  The
graphs $G_1$ and $G_4$ are also $(3,1)$-edge-critical. The next
theorem determines all the $(3,1)$-edge-critical triangle-free
graphs of order $10$.

\begin{theorem}
\label{thm:five} A triangle-free graph $G$ of order 10 is
$(3,1)$-edge-critical if and only if it is isomorphic to $G_5$ or
$G_{1}\cup{K_1}$ or $G_{4}\cup{K_1}$.

\end{theorem}
\smallskip
\noindent{\it Proof.} Let $G$ be a $(3,1)$-edge-critical
triangle-free graph of order10. By Theorem~\ref{thm:four}, either
$G\cong{G_{5}}$ or there is a vertex $u^*$ in $G$ such that
$G-u^*\cong{G_{i}}$ for $1 \le{i}\le{4}$. Clearly the vertex $u^*$
must be an isolated vertex and $i$ must be equal to 1 or 4. Hence
$G$ is isomorphic to $G_5$ or $G_{1}\cup{K_1}$ or $G_{4}\cup{K_1}$.

It is easy to see that $G_{1}\cup{K_1}$ and $G_{4}\cup{K_1}$ are
$(3,1)$-edge-critical. To complete the proof of the theorem we will
show that $\chi_{1}(G_{5}-e)=2$ for every edge $e$ of $G_{5}$.
Clearly $\chi_{1}(G_{5}-e)\ge {2}$ for every edge $e$ of $G_{5}.$

Suppose that $e=(u,u_{1}).$ The sets
$$X_1=\{u,v,u_1,z_1,z_2\}{\mbox{\rm\  and\ }}
V(G_5)-X_1=\{u_2,u_3,u_4,u_5,z\}$$ are 1-independent and hence
provide a (2,1)-colouring of $G_{5}-e$. The edges $(u,u_{2})$,
$(v,u_{1})$ and $(v,u_{2})$ are similar to $(u,u_{1})$ and it is
easy to show that the removal of any one of these edges reduces
$\chi_1(G_5).$

Next suppose that $e=(v,u_{3})$ or $(u,u_{3}).$ The sets
$$X_1=\{u,v,u_3,z\} {\mbox{\rm\  and\ }}
V(G_{5})-X_1=\{u_1,u_2,u_4,u_5,z_1,z_2\}$$ provide a partition of
$V(G_{5}-e)$ into 1-independent sets and hence
$\chi_{1}(G_{5}-e)=2.$ Suppose that $e=(v,u_{4})$ or $(u,u_{4}).$
The sets
$$X_2=\{u,v,u_4,z,z_2\}{\mbox{\rm\  and\ }}
V(G_5)-X_2=\{u_1,u_2,u_3,u_5,z_1\}$$ are 1-independent and hence
$\chi_{1}(G_{5}-e)=2$. Similarly if $e=(v,u_{5})$ or $(u,u_{5})$ the
sets $$X_3=\{u,v,u_5,z,z_1\}{\mbox{\rm\  and\ }}
V(G_5)-X_3=\{u_1,u_2,u_3,u_4,z_2\}$$ are 1-independent and so
$\chi_{1}(G_{5}-e)=2$.

If $e=(u_{3},z_1)$ (resp. $(u_{3},z_2)$), then the sets
$X_1=\{u_1,u_2,u_3,u_4,u_5,z_1$ (resp. $z_2)\}$ and $V(G_{5})-X_1$
provide a $(2,1)$-colouring of $G_{5}-e$. If $e=(u_{4},z_1)$ (resp.
$(u_{5},z_2)$), then the sets $X_2=\{u_1,u_2,u_3,u_4,u_5,z_1$ (resp.
$z_2)\}$ and $V(G_{5})-X_2$ provide a $(2,1)$-colouring of
$G_{5}-e$.

Now if $e=(z,z_i)$ for $i=1$ or $2$ the sets $X_1=\{u,v,z,z_1,z_2\}$
and $V(G_{5})-X_1$ provide a $(2,1)$-colouring of $G_{5}-e$.

Finally if $e=(z,u_i)$ for $i=1$ or $2$ the sets
$$X_1=\{u,v,z_1,z_2\}{\mbox{\rm\  and\ }} V(G_{5})-X_1$$ provide a $(2,1)$-colouring
of $G_{5}-e$.

Thus we have shown that for each $e$ in $G_5$ we have
$\chi_{1}(G_{5}-e)=2$.

This completes the proof of the theorem. $\square$

\medskip

It is easy to see that if a graph with no isolated vertices is
$(3,1)$-edge-critical then it is also $(3,1)$-critical. From
Theorem~\ref{thm:four} it follows that if $G\not\cong{G_5}$ is a
triangle-free graph of order 10 with $\chi_1(G)={3}$ then $G$ is not
$(3,1)$-critical. Hence we have the following theorem.

\begin{theorem}
\label{thm:six} A triangle-free graph $G$ of order 10 is
$(3,1)$-critical if and only if it is isomorphic to $G_5$ given in
Figure~\ref{fig:g5}.

\end{theorem}


\begin{thebibliography}{}
%
%
\bibitem{AAK12b}
Nirmala Achuthan, N.R. Achuthan and G. Keady, On minimal
triangle-free planar graphs with prescribed 1-defective chromatic
number, (to be submitted)

\bibitem{AAS11}
Nirmala Achuthan, N.R. Achuthan, M. Simanihuruk, On minimal
triangle-free graphs with prescribed $k$-defective chromatic number,
{\it Discrete Mathematics} {\bf 311},  1119--1127 (2011).

\bibitem{Avis}
D. Avis, On minimal 5-chromatic triangle-free graphs, {\it J. Graph
Theory}, {\bf 3} , 397--400 (1987).

\bibitem{CLZ}
G. Chartrand, L. Lesniak and P. Zhang, {\it Graphs and Digraphs},
5th ed., Chapman and Hall, 2011.

\bibitem{Chvatal}
V. Chv\'{a}tal, The minimality of the Mycielski graph,
 {\it Graphs and Combinatorics}, Springer-Verlag,
 Berlin, Lecture Notes in Mathematics 406 ,  243--246 (1973).

\bibitem{CGJ}
L. Cowen, W. Goddard, and C.R. Jesurum, Defective coloring
revisited, {\it J. Graph Theory} {\bf 24} , 205--219 (1997).

\bibitem{Fr} 
M. Frick, A survey of (m, k)-colorings, {\it Annals of Discrete
Mathematics} {\bf 55} ,  45--58 (1993).

\bibitem{GH}
J. Gimbel and C. Hartman, Subcolorings and the subchromatic number
of a graph, {\it Discrete Mathematics} {\bf 272} , 139--154 (2003).

\bibitem{HM}
 D. Hanson and G. MacGillivray, On small triangle-free graphs,
{\it Ars Combin.} {\bf 35} , 257--263 (1993).

\bibitem{HS} %
G. Hopkins and W. Staton, Vertex partitions and k-small subsets of
graphs, {\it Ars Combin.} {\bf  22} , 19--24 (1986).

\bibitem{JR}
T. Jensen and G.F. Royle, Small graphs with chromatic number 5: A
computer search, {\it J. Graph Theory} {\bf 19} ,  107--116 (1995).

\bibitem{Lovasz}
L. Lov\`asz, On the compositions of graphs, {\it Studia Sci. Math.
Hungar.} {\bf 1} ,  237--238 (1966).

\bibitem{SAA97b}
 M. Simanihuruk, Nirmala Achuthan, N.R. Achuthan,
On minimal triangle-free graphs with prescribed 1-defective
chromatic number, {\it Australas. J. Combin.} {\bf 16} , 203--227
(1997).

\bibitem{Wo} 
D. Woodall, Improper colourings of graphs, in R. Nelson and R.J.
Wilson eds., {\it Graph Colourings}, Longman Scientific and
Technical (1990).

\end{thebibliography}
\end{document}